\documentclass[12pt]{amsart}
\usepackage{amsmath,amssymb,amsfonts,amsthm}
\usepackage[all,cmtip]{xy}
\usepackage{fullpage}
\DeclareMathAlphabet{\mathpzc}{OT1}{pzc}{m}{it}
\usepackage[francais,english]{babel}
\begin{document}
\theoremstyle{plain}
\newtheorem{thm}{Theorem}[section]
\newtheorem{lem}[thm]{Lemma}
\newtheorem{cor}[thm]{Corollary}
\newtheorem{prop}[thm]{Proposition}
\newtheorem{rem}[thm]{Remark}
\newtheorem{defn}[thm]{Definition}
\newtheorem{ex}[thm]{Example}
\newtheorem{ques}[thm]{Question}
\newtheorem{fact}[thm]{Fact}
\newtheorem{conj}[thm]{Conjecture}
\numberwithin{equation}{subsection}
\def\theequation{\thesection.\arabic{equation}}
\newcommand{\mc}{\mathcal}
\newcommand{\mb}{\mathbb}
\newcommand{\surj}{\twoheadrightarrow}
\newcommand{\inj}{\hookrightarrow}
\newcommand{\zar}{{\rm zar}}
\newcommand{\an}{{\rm an}} 
\newcommand{\red}{{\rm red}}
\newcommand{\codim}{{\rm codim}}
\newcommand{\rank}{{\rm rank}}
\newcommand{\Ker}{{\rm Ker \,}}
\newcommand{\Pic}{{\rm Pic}}
\newcommand{\Div}{{\rm Div}}
\newcommand{\Hom}{{\rm Hom}}
\newcommand{\im}{{\rm im \,}}
\newcommand{\Spec}{{\rm Spec \,}}
\newcommand{\Sing}{{\rm Sing}}
\newcommand{\Char}{{\rm char}}
\newcommand{\Tr}{{\rm Tr}}
\newcommand{\Gal}{{\rm Gal}}
\newcommand{\Min}{{\rm Min \ }}
\newcommand{\Max}{{\rm Max \ }}
\newcommand{\CH}{{\rm CH}}
\newcommand{\pr}{{\rm pr}}
\newcommand{\cl}{{\rm cl}}
\newcommand{\gr}{{\rm Gr }}
\newcommand{\Coker}{{\rm Coker \,}}
\newcommand{\id}{{\rm id}}
\newcommand{\Rep}{{\bold {Rep} \,}}
\newcommand{\Aut}{{\rm Aut}}
\newcommand{\GL}{{\rm GL}}
\newcommand{\Bl}{{\rm Bl}}
\newcommand{\Jab}{{\rm Jab}}
\newcommand{\alb}{\rm Alb}
\newcommand{\NS}{\rm NS}
\newcommand{\sA}{{\mathcal A}}
\newcommand{\sB}{{\mathcal B}}
\newcommand{\sC}{{\mathcal C}}
\newcommand{\sD}{{\mathcal D}}
\newcommand{\sE}{{\mathcal E}}
\newcommand{\sF}{{\mathcal F}}
\newcommand{\sG}{{\mathcal G}}
\newcommand{\sH}{{\mathcal H}}
\newcommand{\sI}{{\mathcal I}}
\newcommand{\sJ}{{\mathcal J}}
\newcommand{\sK}{{\mathcal K}}
\newcommand{\sL}{{\mathcal L}}
\newcommand{\sM}{{\mathcal M}}
\newcommand{\sN}{{\mathcal N}}
\newcommand{\sO}{{\mathcal O}}
\newcommand{\sP}{{\mathcal P}}
\newcommand{\sQ}{{\mathcal Q}}
\newcommand{\sR}{{\mathcal R}}
\newcommand{\sS}{{\mathcal S}}
\newcommand{\sT}{{\mathcal T}}
\newcommand{\sU}{{\mathcal U}}
\newcommand{\sV}{{\mathcal V}}
\newcommand{\sW}{{\mathcal W}}
\newcommand{\sX}{{\mathcal X}}
\newcommand{\sY}{{\mathcal Y}}
\newcommand{\sZ}{{\mathcal Z}}
\newcommand{\A}{{\mathbb A}}
\newcommand{\B}{{\mathbb B}}
\newcommand{\C}{{\mathbb C}}
\newcommand{\D}{{\mathbb D}}
\newcommand{\E}{{\mathbb E}}
\newcommand{\F}{{\mathbb F}}
\newcommand{\G}{{\mathbb G}}
\renewcommand{\H}{{\mathbb H}}
\newcommand{\I}{{\mathbb I}}
\newcommand{\J}{{\mathbb J}}
\newcommand{\M}{{\mathbb M}}
\newcommand{\N}{{\mathbb N}}
\renewcommand{\P}{{\mathbb P}}
\newcommand{\Q}{{\mathbb Q}}
\newcommand{\R}{{\mathbb R}}
\newcommand{\T}{{\mathbb T}}
\newcommand{\V}{{\mathbb V}}
\newcommand{\W}{{\mathbb W}}
\newcommand{\X}{{\mathbb X}}
\newcommand{\Y}{{\mathbb Y}}
\newcommand{\Z}{{\mathbb Z}}
\newcommand{\Nwt}{{\rm Nwt}}
\newcommand{\Hdg}{{\rm Hdg}}
\newcommand{\ind}{{\rm ind \,}}
\newcommand{\Br}{{\rm Br}}
\newcommand{\inv}{{\rm inv}}
\newcommand{\Nm}{{\rm Nm}}
\newcommand{\Griff}{{\rm Griff}}
\newcommand{\Image}{\rm Im \,}
\newcommand{\Ev}{\rm Ev \,}
\title[Stable $\A^1$-connectedness]{Stable $\A^1$-connectedness}
\author{{Nguyen Le Dang Thi}}
\date{27. 08. 2012}          
\begin{abstract}
We prove in this note a stabilized version of a conjecture on $\A^1$-connectedness of A. Asok and B. Dorant \cite[Conj. 2.18]{AD09}. For the stabilized version of this conjecture, we introduce the notion of stable $\A^1$-connectedness, which is can be seen as the stabilization of $\A^1$-connectedness. The notion of stable $\A^1$-connectedness is in connection to the existence of zero cycles of degree one rather than rational points.     
\end{abstract}
\maketitle
\section{Introduction}
Let $k$ be a base field. We denote by $\sH(k)$ the unstable $\A^1$-homotopy category constructed in \cite{MV01}. A $k$-space $\sX$ is $\A^1$-connected, if $\pi_0^{\A^1}(\sX) \stackrel{\cong}{\rightarrow} \pi_0^{\A^1}(\Spec k) = \Spec k$, where we denote by $\pi_0^{\A^1}(\sX)$ the Nisnevich sheaf 
$$\pi_0^{\A^1}(\sX) = a_{Nis}[U \mapsto \Hom_{\sH(k)}(U,\sX)].$$ 
It is conjectured by A. Asok and B. Dorant that:  
\begin{conj}\label{conj1}\cite[Conj. 2.18]{AD09} 
Let $X$ be a smooth $\A^1$-connected scheme and $U \subset X$ an open subscheme such that $\codim(X \setminus U) \geq 2$. If $k$ is an infinite field, then $U$ is also $\A^1$-connected.
\end{conj}
\noindent Remark this conjecture is answered by A. Asok in \cite[Prop. 3.8]{A12}. In this note, we establish the stabilized  version of this conjecture. Denote by $\mathcal{SH}(k)$ resp. $\mathcal{SH}_s(k)$ the $\P^1$-stable resp. the $S^1$-stable homotopy category.  
\begin{defn}\label{def2}
Let $\sX \in Spc(k)$ be a $k$-space. $\sX$ is stable $\A^1$-connected, if the canonical map of sheaves 
$$\pi_0^{st}(\Sigma^{\infty}_{\P^1}\sX_+) \rightarrow \pi_0^{st}(\Sigma^{\infty}_{\P^1}\Spec k_+) $$
is an isomorphism.   
\end{defn}
\noindent In \cite{A12}, A. Asok characterized the notion of unstable $\A^1$-connectedness (cf. \cite[Thm. 4.15]{A12}). A smooth proper $k$-scheme $X$ is $\A^1$-connected if and only if the canonical map of sheaves $\bold{H}_0^{\A^1}(X) \rightarrow \bold{H}^{\A^1}_0(\Spec k) = \Z$ is an isomorphism. Remark that if $X$ is $\A^1$-connected, then $X(k) \neq \emptyset$. One has Hurewicz functors on homotopy categories $\mathcal{SH}_s(k) \rightarrow D_{\A^1}^{eff}(k)$ resp. $\mathcal{SH}(k) \rightarrow D_{\A^1}(k)$, which induce morphisms of sheaves $\pi_i^{st}(\Sigma^{\infty}_s \sX_+) \rightarrow \bold{H}_i^{\A^1}(\sX)$ resp. $\pi_i^{st}(\Sigma^{\infty}_{\P^1} \sX_+) \rightarrow \bold{H}_i^{st\A^1}(\sX)$ for a $k$-space $\sX$. The stable Hurewicz theorem said that (cf. \cite[Thm. 2.3.8]{AH11}) the Hurewicz morphism 
$$\pi_0^{st}(\Sigma^{\infty}_{\P^1} \sX_+) \rightarrow \bold{H}_0^{st\A^1}(\sX)$$     
is an isomorphism convariantly functorial in $\sX$. In particular, if $\sX = \Spec k$, one has $\bold{K}^{MW}_0 \cong \bold{H}_0^{st\A^1}(\bold{1}_k)$. So the characterization for the notion of stable $\A^1$-connectedness is: A $k$-space $\sX$ is stable $\A^1$-connected if and only if $\bold{H}^{st\A^1}_0(\sX) \stackrel{\cong}{\rightarrow} \bold{K}^{MW}_0 = \bold{GW}$, where we denote by $\bold{GW}$ the unramified Grothendieck-Witt sheaf. On the other hand, if $k$ is an infinite perfect field of characteristic unequal $2$, then by \cite[Lem. 4.4.4 \& Thm. 4.4.5]{AH11} the structure map $\bold{H}^{st\A^1}_0(X) \rightarrow \bold{H}^{st\A^1}_0(\Spec k)$ is a split epimorphism if and only $X$ has a zero cycle of degree one, where $X$ is a smooth proper scheme over $k$. So if $X$ is a smooth proper stable $\A^1$-connected $k$-scheme, then the index of $X$ over $k$ is necessary equal to $1$. Of course, the existence of a rational point over $k$ will imply the existence of a zero cycle of degree one over $k$. Now we state the following result: 
\begin{thm}\label{thm3} 
Let $k$ be an infinite perfect field of $char(k) \neq 2$. Let $X$ be a smooth projective stable $\A^1$-connected $k$-scheme. If $U \subset X$ is an open subscheme such that $\codim(X\setminus U) \geq 2$, then $U$ is also stable $\A^1$-connected. 
\end{thm}   
\section{A simple characterization of stable $\A^1$-connectedness}
\noindent We fix an infinite field $k$. Very similar to $\mathcal{SH}(k)$, whose homological $t$-structure is described by F. Morel in \cite{Mor04}, the $t$-structure of $D_{\A^1}(k)$ is given in \cite[Thm. 3.3.3]{AH11}, which is the category of homotopy modules.  Recall (cf. \cite[Defn. 3.3.1]{AH11}) that a homotopy module is a pair $(F_*,\varepsilon_*)$, where $F_*$ is a $\Z$-graded strictly $\A^1$-invariant sheaf and for each $n \in \Z$ one has an isomorphism of abelian sheaves 
$$\varepsilon_n : F_n \rightarrow (F_{n+1})_{-1}.$$
One denotes by $\mathcal{A}b^{st\A^1}(k)$ the category of homotopy modules, whose morphisms are homogeneous natural transformations of $\Z$-graded sheaves compatible with the isomorphisms above. 
One has 
\begin{prop}\label{prop5}\cite[Prop. 3.4.1]{AH11}
Let $F_* \in \mathcal{A}b^{st\A^1}(k)$. One has canonical bijection contravariantly functoral in $\sX$:
$$H^0_{Nis}(\sX,F_0) \stackrel{\cong}{\rightarrow} \Hom_{\mathcal{A}b^{st\A^1}(k)}(\bold{H}_0^{st\A^1}(\sX)_*,F_*). $$ 
\end{prop} 
\noindent This implies immediately the following 
\begin{prop}\label{prop6}
Let $\sX$ be a $k$-space. $\sX$ is stable $\A^1$-connected if and only if for any $F_* \in \mathcal{A}b^{st\A^1}(k)$, one has a bijection 
$$H^0_{Nis}(\sX,F_0) \cong \Hom_{\mathcal{A}b^{st\A^1}}(\bold{GW}_*,F_*) $$
\end{prop}
\begin{ex}\label{ex6}{\rm 
If $X/k$ is a smooth projective $k$-rational variety, then $X$ is stable $\A^1$-connected. Indeed, by \cite[Prop. 3.4.2]{AH11}, the sheaf $\bold{H}_0^{st\A^1}(X)$ is a birational invariance in sense that if $X$ and $X'$ are birationally equivalent smooth proper $k$-varieties, then $\bold{H}_0^{st\A^1}(X) \cong \bold{H}_0^{st\A^1}(X')$. Thus, it is enough to see that $\bold{H}_0^{st\A^1}(\P^n_k) \cong \bold{H}_0^{st\A^1}(\bold{1}_k)$. This follows immediately from the distinguished triangle 
$$C_*^{st\A^1}(\P^n) \rightarrow C_*^{st\A^1}(\P^{n+1}) \rightarrow \Z(n+1)[2n+2] \rightarrow C_*^{st\A^1}(\P^n)[1] $$
}
\end{ex}
\begin{rem}\label{rem7}{\rm 
The same proof as in \cite[Thm. 3.9]{A12}, one can even show that $\bold{H}_0^{st\A^1}(-)$ is a stably $k$-birational invariance of smooth proper $k$-varieties, i.e. if $X$ is stably $k$-birationally equivalent to $X'$, then $\bold{H}_0^{st\A^1}(X) \cong \bold{H}_0^{st\A^1}(X')$. So as in \cite[Thm. 1.4]{N12}, $\bold{H}_0^{st\A^1}(X)$ yields a well-defined invariant on the Grothendieck ring $K_0(Var/k)/\bold{L}$, where $\bold{L}$ denotes the Lefschetz class $[\A^1_k]$.    
}
\end{rem}
\section{Proof of \ref{thm3}}
\noindent Throughout this section, we fix a base field $k$, which is an infinite perfect field of $char(k) \neq 2$. Let $X$ be a smooth projective stable $\A^1$-connected $k$-scheme. We write $Z = X \setminus U$ for the complement and $c = \codim_X(Z) \geq 2$. Assume first of all $Z$ is smooth, remark that we can use the stable $\A^1$-connectivity theorem (cf. \cite{Mor05}) of F. Morel to prove immediately the theorem, however we give here in the first part another proof including duality and calculation with $\A^1$-cohomology. Consider then the Gysin triangle in $D_{\A^1}(k)$ (see \cite{Deg08})  
$$C_*^{st\A^1}(U) \rightarrow C_*^{st\A^1}(X) \rightarrow C_*^{st\A^1}(Z)(c)[2c] \rightarrow C_*^{st\A^1}(U)[1].$$
This yields an exact sequence of abelian groups 
\begin{multline*} 
\cdots \rightarrow \Hom_{D_{\A^1}(k)}(C_*^{st\A^1}(\Spec L), C_*^{st\A^1}(Z)(c)[2c-1]) \rightarrow \Hom_{D_{\A^1}(k)}(C_*^{st\A^1}(\Spec L),C_*^{st\A^1}(U)) \\ \rightarrow \Hom_{D_{\A^1}(k)}(C_*^{st\A^1}(\Spec L),C_*^{st\A^1}(X)) \rightarrow \Hom_{D_{\A^1}(k)}(C_*^{st\A^1}(\Spec L),C_*^{st\A^1}(Z)(c)[2c]) \rightarrow \cdots 
\end{multline*}
where $L \in \sF_k$ is a separable finitely generated field extension over $k$. By definition, we have $\Hom_{D_{\A^1}(k)}(C_*^{st\A^1}(\Spec L), C_*^{st\A^1}(X)) = \bold{H}_0^{st\A^1}(X)(L)$ and analogously for $U$. Now we prove that the other two terms in the exact sequence above are trivial. By Atiyah-Spanier-Whitehead duality in $D_{\A^1}(k)$ (see \cite[\S 3.5]{AH11}), we can write 
$$\Hom_{D_{\A^1}(k)}(C_*^{st\A^1}(\Spec L),C_*^{st\A^1}(Z)(c)[2c-1]) \cong H^{2(c+n_Z+d_Z)-1}_{st\A^1}(Th(V_{Z_L}),\Z(c+n_Z+d_Z)),$$
where the vector bundle $V_Z$ comes from \cite[Thm. 2.11]{Voe03} and $n_Z$ denotes its rank. Here we write $d_Z = \dim(Z)$ and $H^{*}_{st\A^1}(-,\Z(*))$ is the $\A^1$-stable cohomology. We compute first of all the unstable $\A^1$-cohomology $H^{2(c+n_Z+d_Z)-1}_{\A^1}(Th(V_{Z_L}),\Z(c+n_Z+d_Z))$. By \cite[Prop. 3.2.5]{AH11}, it is nothing but the hypercohomology $\H^*_{Nis}(-,\Z_{\A^1}(*))$. Consider the hypercohomology spectral sequence 
$$E^{p,q}_2 = H^p_{Nis}(-,\underline{H}^q(\Z_{\A^1}(*))) \Rightarrow \H^{p+q}_{Nis}(-,\Z_{\A^1}(*)).$$
As $\underline{H}^q(\Z_{\A^1}(*)) = 0$ for $q > *$ (a consequence of Morel's stable $\A^1$-connectivity theorem) and $\bold{K}_*^{KW} = \underline{H}^*(\Z_{\A^1}(*))$ for $* > 0$ (a consequence of stable Hurewicz theorem), we have a canonical homomorphism 
\begin{equation}\label{eq4}
H_{\A^1}^{2(c+n_Z+d_Z)}(Th(V_{Z_L}),\Z(c+n_Z+d_Z)) \rightarrow H^{c+d_Z+n_Z-1}_{Nis}(Th(V_{Z_L}),\bold{K}^{MW}_{c+n_Z+d_Z}). 
\end{equation} 
Since $c \geq 2$, which implies that $2(c+n_Z+d_Z)-1 > cd_{Nis}(Th(V_{Z_L})) + c + n_Z + d_Z$ (remark that $cd_{Nis}(Th(V_{Z_L})) \leq n_Z + d_Z$, cf. \cite[Lem. 3.1.7]{AH11}), so by \cite[Cor. 3.2.6]{AH11} the above homomorphism is an isomorphism. But the group in the right hand side of the homomorphism \ref{eq4} is trivial because of cohomological dimension reason. So we see that the unstable $\A^1$-cohomology $H_{\A^1}^{2(c+n_Z+d_Z)}(Th(V_{Z_L}),\Z(c+n_Z+d_Z))$ is trivial. Let $j>0$ be a natural number, we consider the groups 
$$H^{2(c+n_Z+d_Z)+j-1}_{\A^1}(Th(V_{Z_L}) \wedge \G_m^{\wedge j},\Z(c+n_Z+d_Z+j)).$$      
These groups can be rewritten as $H^{2(c+n_Z+d_Z)+j-1}_{\A^1}(Th(V_{Z_L}),\Z(c+n_Z+d_Z+j))_{-(j)}$. So the same argument with the spectral sequence of hypercohomology as above yields 
$$H^{2(c+n_Z+d_Z)+j-1}_{\A^1}(Th(V_{Z_L}),\Z(c+n_Z+d_Z+j)) \cong H^{c+n_Z+d_Z-1}_{Nis}(Th(V_{Z_L}),\bold{K}^{MW}_{c+n_Z+d_Z+j}),$$
which are trivial again by cohomological dimension reason ($c \geq 2$). So we may conclude now that the stable $\A^1$-cohomology $ H^{2(c+n_Z+d_Z)-1}_{st\A^1}(Th(V_{Z_L}),\Z(c+n_Z+d_Z))$ is trivial. Similarly, the group $\Hom_{D_{\A^1}(k)}(C_*^{st\A^1}(\Spec L),C_*^{st\A^1}(Z)(c)[2c])$ must also vanish. We see then $\bold{H}_0^{st\A^1}(U)(L) \cong \bold{H}_0^{st\A^1}(X)(L) \cong GW(L)$. In general, if $Z$ is not smooth, since $k$ is perfect, there is a finite stratification of closed subschemes of $Z$ 
$$\emptyset \subset Z_{d_X} \subset Z_{d_X-1} \subset \cdots \subset Z_{c+1} \subset Z_c = Z,$$    
such that each $Z_i \setminus Z_{i+1}$ is smooth in $X \setminus Z_{i+1}$. So we may prove recursively and come back to the case when $Z = Z_{d_X}$, which follows from the previous case, since $0$-dimensional schemes are smooth. More precisely, let us consider the Gysin triangle 
$$C_*^{st\A^1}(X-Z_i) \rightarrow C_*^{st\A^1}(X-Z_{i+1}) \rightarrow C_*^{st\A^1}(Z_i - Z_{i+1})(c_i)[2c_i], $$
where all $c_i \geq 2$. As before, we apply the functor $\Hom_{D_{\A^1}(k)}(C_*^{st\A^1}(\Spec L),-)$ and consider then the group $\Hom_{D_{\A^1}(k)}(C_*^{st\A^1}(\Spec L),C_*^{st\A^1}(Z_i - Z_{i+1})(c_i)[2c_i - j])$  for $j= 0,1$. These groups can be rewritten as 
\begin{multline*} 
\Hom_{D_{\A^1}(k)}(C_*^{st\A^1}(\Spec L)[j-c_i],C_*^{st\A^1}(Z_i - Z_{i+1}\wedge \G_m^{\wedge c_i}) ) = \\ = colim_n \Hom_{D_{\A^1}^{eff}(k)}(\tilde{C}_*^{\A^1}(\Spec L_+ \wedge \G_m^{n})[j-c_i],\tilde{C}_*^{\A^1}((Z_i-Z_{i+1})_+ \wedge \G_m^{\wedge n+c_i})) = \\ = colim_n H_{j-c_i}^{\A^1}((Z_i-Z_{i+1}) \wedge \G_m^{\wedge c_i+n})_{-(n)}(\Spec L). 
\end{multline*}
The last group must vanishing according to the stable $\A^1$-connectivity theorem of F. Morel \cite{Mor05}, as $c_i \geq 2$ and $j = 0,1$. 
\bibliographystyle{plain}
\renewcommand\refname{References}

\end{document}